\documentclass[12pt]{amsart}
\usepackage{amssymb}
\usepackage[all]{xy}

\title{Higher Nash blowups}
\author{Takehiko Yasuda}
\address{Research Institute for Mathematical Sciences, Kyoto University, Kyoto, 606-8502, Japan}
\date{\today}
\email{takehiko@kurims.kyoto-u.ac.jp}

\theoremstyle{plain}
\newtheorem{thm}{Theorem}[section]
\newtheorem{conj}[thm]{Conjecture}
\newtheorem{prop}[thm]{Proposition}
\newtheorem{cor}[thm]{Corollary}
\newtheorem{lem}[thm]{Lemma}

\theoremstyle{definition}
\newtheorem{defn}[thm]{Definition}
\newtheorem{expl}[thm]{Example}

\theoremstyle{remark}
\newtheorem{rem}[thm]{Remark}

\def\AA{\mathbb A}
\newcommand{\CC}{\mathbb C}

\newcommand{\xx}{\mathbf x}
\newcommand{\yy}{\mathbf y}

\newcommand{\bs}{\mathbf{s}}
\newcommand{\bu}{\mathbf{u}}
\newcommand{\ba}{\mathbf{a}}

\newcommand{\bt}{\mathbf{t}}
\newcommand{\NN}{\mathbb N}
\newcommand{\No}{\mathbb{N}_0}

\newcommand{\cI}{\mathcal{I}}

\newcommand{\cK}{\mathcal{K}}

\newcommand{\cM}{\mathcal{M}}
\newcommand{\cN}{\mathcal{N}}
\newcommand{\cO}{\mathcal{O}}
\newcommand{\cP}{\mathcal{P}}

\newcommand{\cZ}{\mathcal{Z}}
\newcommand{\cW}{\mathcal{W}}

\newcommand{\fa}{\mathfrak{a}}
\newcommand{\fb}{\mathfrak{b}}
\newcommand{\fc}{\mathfrak{c}}
\newcommand{\fj}{\mathfrak{j}}

\newcommand{\fm}{\mathfrak{m}}

\newcommand{\ann}{\mathrm{ann}}

\newcommand{\length}{\mathrm{length}\,}
\newcommand{\ord}{\mathrm{ord}\,}

\newcommand{\Spec}{\mathrm{Spec}\,}

\newcommand{\pr}{\mathrm{pr}}

\newcommand{\sm}{\mathrm{sm}}

\newcommand{\Hilb}{\mathbf{Hilb}}
\newcommand{\Grass}{\mathbf{Grass}}

\newcommand{\Nash}{\mathbf{Nash}}

\newcommand{\red}{\mathrm{red}}

\allowdisplaybreaks[2]

\begin{document}

\begin{abstract}
For each non-negative integer $n$,
we define the $n$-th Nash blowup of an algebraic variety, and
call them all higher Nash blowups.
When $n=1$, it coincides with the classical Nash blowup.
We study higher Nash blowups of curves in detail
and prove that any curve in characteristic zero
can be desingularized by its $n$-th Nash blowup 
with $n$ large enough.
\end{abstract}

\maketitle

\section*{Introduction}

The classical Nash blowup of an algebraic variety is the parameter space of
the tangent spaces of smooth points and their limits,
and the normalized Nash blowup is the Nash blowup followed by the normalization.
It is natural to ask whether the iteration of
Nash blowups or normalized Nash blowups leads to a smooth variety.
There are works on this question, by Nobile \cite{Nobile}, Rebassoo \cite{Rebassoo},
 Gonz\'alez-Sprinberg \cite{Gonzalez-Sprinberg}, Hironaka \cite{Hironaka-Nash} and Spivakovsky \cite{Spivakovsky}. If the answer is affirmative, we obtain a canonical way to resolve singularities.

In this paper, we make a similar but different approach to a resolution of singularities.
Let $X$ be an algebraic variety over an algebraically closed field $k$.
For a point $x \in X$, we denote by $x^{(n)}$ its $n$-th infinitesimal neighborhood,
that is, if $(\cO_{X,x} ,\fm_x)$ is the local ring at $x$, 
the closed subscheme $\Spec \cO_{X,x} /\fm_x^{n+1} \subseteq X$. 
If $x$ is a smooth point, being an Artinian subscheme of length $ \binom{n+d}{d} $,
$x^{(n)}$ corresponds to a point $[x^{(n)}]$ of the Hilbert scheme $\Hilb _{\binom{n+d}{d}} (X)$
of $\binom{n+d}{d}$ points of $X$.
We define the {\em $n$-th Nash blowup} of $X$, denoted $\Nash _n (X)$, 
to be the closure of the set 
\[
\{ (x,[x ^{(n)}]) | x \text{ smooth point of }X\}
\]
 in $X \times _{k} \Hilb _{\binom{n+d}{d}} (X)$.
We also call it a {\em higher Nash blowup} of $X$.
The first projection restricted to $\Nash _{n} (X)$
\[
 \pi_{n} :\Nash _{n} (X) \to X
\]
is a projective birational morphism which is an isomorphism
over the smooth locus of $X$.
The first Nash blowup is canonically isomorphic to the classical Nash blowup (see Proposition \ref{prop-isoms-coherent-Nash}).
Every point of $\Nash _ n (X)$ corresponds to an Artinian subscheme $Z$ of $X$
which is set-theoretically a single point.

If $\Nash _{n}'(X)$ is the closure of  $\{ [x ^{(n)}] | x \text{ smooth point of }X \}$
in $\Hilb _{\binom{n+d}{d}} (X)$, then there exists a natural morphism
$ \Nash _{n}(X) \to \Nash _{n} (X')$, $(x,[Z]) \mapsto [Z]$.
Thus $\Nash _n (X)$ is identified with the set of
 the $n$-th infinitesimal neighborhoods of
smooth points and their limits. 
We can also construct higher Nash blowups by using the relative Hilbert scheme
or the Grassmaniann schemes of coherent sheaves. 
The last construction
is essentially the same as  a special case of Oneto and Zatini's Nash blowup
associated to a coherent sheaf \cite{Oneto-Zatini-Nash}.

The problem that interests us is of course
{\em whether varieties can be desingularized by higher Nash blowups.}
Concerning this problem, we formulate a conjecture as follows:

\begin{conj}\label{Yasuda-conjecture}
Suppose that $k$ has characteristic zero.
Let $X$ be a variety of dimension $d$, 
$J^{(d-1)}$ the $(d-1)$-th neighborhood of the Jacobian subscheme  $J \subseteq X$
(that is, the closed subscheme defined by the $d$-th power of the Jacobian ideal sheaf).
Let $[Z] \in \Nash_n (X)$ with $Z \nsubseteq J^{(d-1)}$.
Then $\Nash _n (X)$ is smooth at $[Z]$. 
\end{conj}

For any closed subscheme $Y \subseteq X$ of dimension $< d$, there exists $n_0$
such that for every $n \ge n_0$  and for every $[Z] \in \Nash _n (X)$, $Z \nsubseteq Y$ 
(Proposition \ref{prop-Hilbert-func-inclusion}). 
Then the conjecture especially says that
$\Nash _n (X)$ are smooth for  $n \gg 0$. 
If the conjecture is true, we obtain a canonical way to resolve singularities by
{\em one step}.

Our first step toward proving the conjecture is a separation of analytic branches.
Let $\hat X := \Spec \hat \cO_{X,x}$ be the completion of a variety $X$ at $x \in X$,
and $\hat X_i$, $i=1,\dots,l$, its irreducible components.
Then we can define higher Nash blowups of $\hat X$ and $\hat X_{i}$, and obtain
\[
 \Nash _n ( X) \times _X \hat X \cong \Nash _n (\hat X) = \bigcup _{i=1}^l \Nash _n (\hat X_i).
\]
Let $\nu : \tilde X \to X$ be the normalization.
The {\em conductor ideal sheaf} is the annihilator ideal sheaf of
the coherent sheaf $\nu _{* } \cO_{\tilde {X}} / \cO_{X}$.
The {\em conductor subscheme} $C \subseteq X$ is the closed subscheme
defined by the conductor ideal sheaf.

\begin{prop}[=Proposition \ref{prop-general-separation}]\label{prop-hypersurface-separation-introduction}
Let $[Z] \in \Nash_{n} (X)$ with $Z \nsubseteq C$.
Then $Z$ is contained in a unique analytic branch of $X$. 
\end{prop}

If $x \in X$ is the support of $Z$ and $\hat X_{i}$ are as above,
then the proposition says that $\Nash _{n} (\hat X_{i})$
are disjoint around $[Z]$. Therefore the study of $\Nash _{n} (X)$
is reduced to that of $\Nash _{n} (\hat X_{i})$.

We study the case of curves in more detail. 
Let $R$ be a local complete Noetherian domain of 
dimension $1$ with coeffiecient field $k$
and $X := \Spec R$. The integral closure of $R$
is (isomorphic to) $k[[x]]$. Then we define a numerical monoid
$S := \{ i | \exists f \in R, \ \ord f =i \}$.
In characteristic 0, we can completely determine when $\Nash _{n} (X)$ is regular\footnote{Since
$\Nash _{n } (X)$ is not of finite type over $k$, we use the term ``regular'' instead of ``smooth''.}
in terms of $S$.

\begin{thm}[=Theorem \ref{thm-curve1}]\label{thm-curve1-intro}
Let $X$ and $S$ be as above. Suppose that $k$ has characteristic zero.
Then
$\Nash _n (X)$ is regular if and only if $ s_n -1 \in S $.
\end{thm}

As a corollary, we prove the following, which implies Conjecture
\ref{Yasuda-conjecture} in dimesion $1$.

\begin{cor}[=Corollary \ref{cor-general-curve}]
Let $X$ be a variety of dimension $1$ over $k$, $C$ its conductor subscheme and 
$[Z] \in \Nash _{n} (X)$. Suppose that $k$ has characteristic 0
and that $Z \nsubseteq C$. Then $\Nash _{n} (X)$ is smooth at
$[Z]$.
\end{cor}

In contrast to the iteration of classical Nash blowups,
each higher Nash blowup is directly constructed from the given variety.
There is no direct relation between $\Nash _{n+1} (X)$ and $\Nash _n (X)$.
In fact, from Theorem \ref{thm-curve1-intro},
we see that even if $\Nash _n (X)$ is regular, $\Nash _{n+1} (X)$ is not generally regular.
So there is no birational morphism $\Nash_ {n+1} (X) \to \Nash _n (X) $
(See Example \ref{expl-5,7}).

\begin{rem}
\begin{enumerate}
\item
There are few evidences for Conjecture \ref{Yasuda-conjecture} in higher dimension.
It is, maybe, safer to replace $J^{(d-1)}$ with $J^{(a_d)}$, where $a_d$ is a positive integer depending only
on $d$. The conjecture is based on the idea that Artinian subschemes protruding much from
the singular locus behave well. 
A similar idea for jets appears in the theory of motivic integration for
singular varieties (see \cite{Denef-Loeser}). 
\item
The conjecture fails, if we replace $J^{(d-1)}$ with $J=J^{(0)}$:
Let 
\[
 X := (x^2+y^2 + z^{n+1} =0) \subseteq \AA^3_\CC
\]
be a surface with an $A_n$-singularity. 
Its Jacobian ideal is $(x,y,z^n) \subseteq \CC[x,y,z]/(x^2+y^2+z^{n+1})$.
Let $A \subseteq X$ be the subscheme defined by the Jacobian ideal, which is isomorphic to
$\Spec \CC[z]/(z^n)$. 
For any $[Z]\in \Nash _1 (X)$, $Z \cong \Spec \CC[s,t]/(s,t)^2$,
and $Z \nsubseteq A$. However the classical Nash blowup of $X$ 
is not generally smooth (see \cite[\S 5.2]{Gonzalez-Sprinberg}).
\item
The conjecture fails also in positive characteristic at least in dimension 1.
Let $X$ be an analytically irreducible curve in characteristic $p >0$.
Then $\Nash _{p^e-1} (X) \cong X$ for  $e \gg 0$ 
(Proposition \ref{prop-positive-char-curve}). 
If $k$ is of characteristic either 2 or 3, and if $X= \Spec k[[x^2,x^3]]$,
 then $\Nash _n (X) \cong X$ for every $n$ (Proposition \ref{prop-cusp-char2}).
 \end{enumerate}
\end{rem}

Nakamura's $G$-Hilbert scheme is also a kind of blowup constructed by using a Hilbert scheme of points.
For an algebraic variety $M$ with an effective action of a finite group $G$,
its $G$-Hilbert scheme $G \text{-} \Hilb (M)$ parameterizes the free orbits and their 
limits in the Hilbert scheme
of points of $M$, and there exists a projective birational morphism $G \text{-} \Hilb (M) \to M/G $. 
Replacing free obits with their $n$-th infinitesimal neighborhoods,
we can define a higher version of $G$-Hilbert scheme, although the author does not know
whether it is interesting.

We can easily generalize the higher Nash blowup to generically smooth morphisms, that is, to
the relative setting, and even more generally to  foliations. The latter was actually what the author
first thought of. 
In sum, the common idea is the following:
Given a variety or a variety with some additional structure (such as a morphism or a foliation), 
the space of some objects uniquely associated to smooth points and their limits is, if well-defined as
a variety, then a modification of the given variety. 
 If the modification is a resolution of singularities, then it should be
  an advantage that the resulting variety is a moduli space of some objects on
  the given variety.  

In Section \ref{section-definition}, we give the definition of higher Nash blowup and
several alternative constructions.
In Section \ref{section-general-properties}, we prove basic properties of higher Nash blowups.
In the final section, we study the case of curves. 

\subsection*{Acknowledgments}

The author would like to thank Shigefumi Mori, Shigeru Mukai, Hisanori Ohashi, Shunsuke Takagi,
Masataka Tomari and Kei-ichi Watanabe
for useful comments.

\subsection*{Conventions}

We work in the category of schemes over an algebraically closed field $k$.
A {\em  point} means a $k$-point.
A {\em variety} means an integral separated scheme of finite type over $k$.
For a closed subscheme $Z \subseteq X$ defined by an ideal $\cI \subseteq \cO_X$,
we denote by $Z^{(n)}$ its {\em $n$-th infinitesimal neighborhoods}, that is,
the closed subscheme defined by $\cI^{n+1}$. 
We denote by $\NN $  the set $\{1,2,\dots\}$ of positive integers 
and by $\No$ the set $\{0,1,2,\dots\}$ of non-negative integers.

\section{Definition and several constructions}\label{section-definition}

\subsection{Definition}

Let $X$ be a variety  of dimension $d$, and $x \in X$ and 
$ x^{(n)} := \Spec \cO_{X,x}/\fm_x ^{n+1}$ its $n$-th infinitesimal neighborhood.
If $X$ is smooth at $x$, then $x^{(n)}$ is an Artinian subscheme of $X$
of length $\binom{d + n }{ n}$. 
Therefore it corresponds to a point
\[
[x^{(n)}] \in 
 \Hilb _{\binom{d+n}{n}}(X)  ,
\]
where $\Hilb _{\binom{d+n}{n}}(X)$ is the Hilbert scheme of 
$\binom{d+n}{n}$ points of $X$.
If $X_\sm$ denotes the smooth locus of $X$, then  we have a map
\begin{align*}
\sigma _n:  X_\sm \to \Hilb_{\binom{d+n}{n}}(X), \ x \mapsto  [x^{(n)}] .
\end{align*}

\begin{lem}\label{lem-morphisms-schemes}
$\sigma _n$ is a morphism of schemes.
\end{lem}

\begin{proof}
Let   $\Delta \subseteq X _\sm \times _k X_\sm $ be the diagonal.
Consider a diagram of the projections restricted to its $n$-th infinitesimal neighborhood
$\Delta^{(n)}$,
\[\xymatrix{
\Delta ^{(n)} \ar[r] ^{\pr _2} \ar[d]_{\pr_1} & X _\sm \\
X_\sm & .
}
\]
For $x \in X_\sm $, 
\[
\pr_2 (\pr_1^{-1}(x)) = x^{(n)}.
\]
Therefore by the definition of Hilbert scheme,
there exists a morphism
\[
 X_\sm \to \Hilb _{\binom{d+n}{n}}(X)
\]
corresponding to the diagram above.
It is identical to $\sigma_n$.
\end{proof}

The graph $\Gamma _{\sigma _{n}} \subseteq X_{\sm} \times _{k} \Hilb_{\binom{d+n}{n}}(X)$
of $\sigma_{n}$ is canonically isomorphic to $X_{\sm}$.

\begin{defn}
We define {\em the $n$-th Nash blowup  of $X$},
denoted $\Nash_n(X)$, to be the closure of $\Gamma_{\sigma_{n}}$ with reduced
scheme structure in $X \times _{k} \Hilb_{\binom{d+n}{n}}(X)$.
\end{defn}

The first projection restricted $\Nash_{n}(X)$,
\[
 \pi_{n} : \Nash _{n} (X) \to X,
\]
is projective and birational.
Moreover it is an isomorphism over $X_{\sm}$.

Let $\Nash _{n}' (X)$ be the closure of 
$\sigma _{n}(X_{\sm})$ in $\Hilb_{\binom{d+n}{n}}(X)$.
Then the second projection $X \times _{k} \Hilb_{\binom{d+n}{n}}(X)$
induces a morphism 
\[
\psi_{n}: \Nash _{n} (X) \to \Nash'_{n}(X) .
\]
This bijectively sends $(x,[Z])$ to $[Z]$.
Thus $\Nash _{n} (X)$ is set-theoretically identified with 
$\Nash'_{n} (X)$, the set of the $n$-th infinitesimal neighborhoods
of smooth points and their limits.
Hereafter we abbreviate $(x,[Z]) \in \Nash _{n} (X)$ as
$[Z] \in \Nash _{n} (X)$.

\subsection{$\psi_{n}$ is an isomorphism in  char.\ 0}

Let $S^{m} X$ denote the $m$-th symmetric product of $X$.
The Hilbert-Chow morphism of \cite{Fogarty} is 
a morphism 
\[
(\Hilb _{m} (X) )_{\red} \to S^{m} X
\]
which assign a closed subscheme $Z \subseteq X$
the associated $0$-cycle.

In characteristic zero, 
 $X$ is embedded into $S^{m} X$ as the small diagonal,
$\{ (x , \dots,x) | x \in X \} \subseteq S^{m} X$.
(In positive characteristic, the diagonal morphism 
$X \to S^{m} X$ is not generally a closed embedding.)
When $m=\binom{d+n}{n}$,  the Hilbert-Chow morphism
restricted to $\Nash '_{n} (X)$,
\[
\pi'_{n} : \Nash _{n} '(X) \to X,
\]
 is a morphism onto $X$.

\begin{prop}
Suppose that $k$ has characteristic 0. 
Then  $\psi_{n}$ is an isomorphism and $\pi_{n} = \pi'_{n} \circ \psi$.
\end{prop}

\begin{proof}
The graph $ \Gamma _{\pi'_{n}}\subseteq X \times _{k} \Nash _{n}' (X) $
of $\pi_{n}'$ is identical to $\Nash _{n} (X)$.
Therefore $\psi_{n}$ is an isomorphism. Now the equality
$\pi_{n} = \pi'_{n} \circ \psi$ is obvious.
\end{proof}

\begin{rem}
In positive characteristic, $\psi_{n}$ is not generally an isomorphism.
For instance, let $X:=\Spec k[x]$. Since $X$ is smooth, 
$\phi_{n}$ is isomorphic to 
$\sigma _{n}:X \to \Nash '_{n}(X) \subseteq \Hilb _{\binom{n+d}{d}} (X)$.
Suppose that $k$ has characteristic $p >0$ and that $p$ divides $n +1$.
Then $\Delta ^{(n)} \times _{X} \Spec k[x]/(x^{2})$ is a trivial 
embedded deformation over $\Spec k[x]/(x^{2})$.
So the corresponding morphism $\Spec k[x]/(x^{2}) \to X \to \Nash '_{n}(X)$
factors as $\Spec k[x]/(x^{2}) \to \Spec k \to \Nash'_{n} (X)$.
It follows that $\phi_{n}$ is not an isomorphism.
\end{rem}

\subsection{Construction with the relative Hilbert scheme}

We can construct higher Nash blowups also by using the {\em relative} Hilbert scheme.
Let $X$ be a variety and $\Delta^{(n)} \subseteq X \times _{k} X$
the $n$-th infinitesimal neighborhood of the diagonal.
Then the restricted first projection
\[
 \pr_{1} :\Delta^{(n)} \to X
\]
is a finite morphism. Its relative Hilbert scheme
\[
\Hilb _{\binom{d + n }{ n}} ( \pr _1 :\Delta ^{(n)} \to X)
\]
for  a constant Hilbert polynomial $\binom{d + n }{ n}$
is a projective $X$-scheme.
It is easy to see that 
\[
 \Hilb _{\binom{d + n }{ n}} ( \pr _1 :\Delta ^{(n)} \to X) \times _{X} X_{\sm} \cong X_{\sm}.
\]

\begin{prop}\label{prop-alternative-definition-of-HB}
The irreducible component of $\Hilb _{\binom{d + n }{ n}} ( \pr _1 :\Delta ^{(n)} \to X)$
dominating $X$ is canonically isomorphic to $\Nash _{n} (X)$.
\end{prop}

\begin{proof}
A closed embedding $\Delta^{(n)} \hookrightarrow X \times _k X$
induces a closed embedding 
\[
 \Hilb _{\binom{n+d}{d}}(\pr _1 : \Delta ^{(n)} \to X) \hookrightarrow \Hilb_{\binom{n+d}{d}} (\pr _1 : X \times _k X\to X).
\]
We also have a closed embedding
\begin{align*}
\Nash_n (X)   \hookrightarrow  X \times_k \Hilb_{\binom{n+d}{d}} (X) = \Hilb_{\binom{n+d}{d}} (\pr _1 : X \times _k X\to X) .
\end{align*}
Then $\Nash _n (X)$ and the irreducible component of 
$\Hilb _{\binom{n+d}{d}}(\pr _1 : \Delta ^{(n)} \to X)  $ dominating $X$
determines the same closed subscheme of $\Hilb_{\binom{n+d}{d}} (\pr _1 : X \times _k X\to X)$.
This proves the assertion.
\end{proof}

\begin{cor}\label{cor-contained-in-nth-nbhd}
Let $[Z] \in \Nash _n (X)$ such that the support of $Z$ is $x$, (that is $\pi_{n}([Z])=x$).
Then $Z \subseteq x ^{(n)}$. 
\end{cor}

\begin{proof}
The subscheme $Z \subseteq X$ is contained in the fiber of $\pr_1:\Delta^{(n)} \to X$ over $x$,
which is exactly $x^{(n)}$.
\end{proof}

\subsection{The Nash blowup associated to a coherent sheaf}

Let $X$ be a reduced 
Noetherian scheme, $ \cM $ a coherent $\cO_X$-module locally free of constant rank $r$
on an open dense subscheme $U \subseteq X$, and
$\Grass _r (\cM)$ the Grassmaniann  of $\cM$ of rank $r$, which is a projective $X$-scheme. 
Then the fiber product $\Grass _r (\cM) \times _X U $ is isomorphic to $U$ by the projection.

\begin{defn}
The closure of $\Grass _r (\cM) \times _X U $ is called the {\em Nash blowup of $X$ associated to} $\cM$
and denoted $\Nash (X,\cM)$ (see \cite{Oneto-Zatini-Nash}). 
\end{defn}

Then the natural morphism $\pi _\cM : \Nash (X,\cM) \to X$ is projective and birational.
When $X$ is a variety and $\cM = \Omega_{X/k}$, then $\Nash (X,\Omega_{X/k})$
is  the classical Nash blowup of $X$.

If $tors \subseteq \pi_\cM ^ * \cM$ denotes the torsion part, then by definition,
$(\pi_\cM ^ * \cM )/ tors$ is locally free. Moreover $\Nash (X,\cM)$ has the following universal property:
If $f:Y \to X$ is a modification with $(f^* \cM)/tors$ locally free, then 
there exists a unique morphism $g : Y \to \Nash (X,\cM)$ with $\pi_\cM \circ g =f $.

Let $\cI _\Delta  \subseteq \cO _{X\times _k X}$ be the ideal sheaf defining the diagonal 
$\Delta \subseteq X \times _k X$.
Put $ \cP ^n _X :=  \cO _{X\times _k X} / \cI _\Delta ^{n +1} $ and 
$\cP_{X,+} ^n := \cI_\Delta / \cI _\Delta^{n+1}$, $ n \in \NN$. The $\cP ^n_X$
is the structure sheaf of $\Delta^{(n)}$ and called 
the {\em sheaf of principal parts of order $n$} of $X$ (see \cite[Def.\ 16.3.1]{EGA}).
 We regard $ \cP ^n _X$ and $\cP_{X,+}^n$ as $\cO_X$-modules through the first projection.
When $X$ is a variety, these are  coherent sheaves.

\begin{prop}\label{prop-isoms-coherent-Nash}
For every variety $X$ and every $n \in \No$, we have canonical isomorphisms,
\[
\Nash _n (X) \cong \Nash (X,\cP^n_X) \cong \Nash (X , \cP^n _{X,+}).
\]
In particular, $\Nash_1 (X)$ is canonically isomorphic to the classical Nash blowup of $X$.
\end{prop}

\begin{proof}
Because of the universal property, if $\cN$ is locally free, then we have a canonical 
isomorphism $ \Nash (X,\cM \oplus \cN) \cong \Nash (X,\cM)$.
In particular, since $\cP^n_X \cong \cO_X \oplus \cP^n _{X,+}$, we have
$ \Nash (X,\cP^n_X) \cong \Nash (X , \cP^n _{X,+})$.

The moduli schemes $\Hilb _{\binom{d+n}{n}}(\pr_{1}:\Delta^{(n)}\to X)$
and $\Grass _{\binom{d+n}{n}} (\cP^{n}_{X})$ represent
equivalent functors. Hence they are canonically isomophic.
It follows that $\Nash _n (X) \cong \Nash (X,\cP^n_X)$.
\end{proof}

\begin{cor}
Let $X$ be a variety of dimension $d$, $n \in \No$, and $r := \binom{n+d}{d}$.
\begin{enumerate}
\item $ \Nash _n (X) \cong \Nash (X, \bigwedge ^{r} \cP^n_X)$.
\item Let $\cK(X)$ be the constant sheaf of rational functions. Fix an isomorphism 
$\bigwedge ^{r} \cP^n_X \otimes _{\cO_X} \cK(X) \to \cK (X)$ and define a homomorphism
\[
 \psi : \bigwedge ^{r} \cP^n_X \to  \bigwedge ^{r} \cP^n_X \otimes _{\cO_X} \cK(X)\to  \cK (X) .
\]
Then  $\Nash _n (X)$ is isomorphic to the blowup of $X$ with respect to a fractional ideal 
$\psi (\bigwedge ^{r} \cP^n_X)$.
\end{enumerate}
\end{cor}

\begin{proof}
These are results due to Oneto and Zatini \cite{Oneto-Zatini-Nash} restricted to the case
where $\cM = \cP^n_X$.
\end{proof}

\subsection{Formal completion}

For a complete local Noetherian ring $S$ with coefficient field $k$,
the module $\Omega_{S/k}$ of K\"ahler differentials is not generally 
finitely generated over $S$, while its completion $\hat \Omega_{S/k}$
is. The latter is usually the suitable one to handle. 
We show those similar facts on its higher version $\hat P_{S}^n$
that are required in  applications to higher Nash blowups.

Let $k[\xx] := k[x_1,\dots ,x_r]$ be a polynomial ring with $r$ variables 
and $R = k[\xx] /\fa$ its quotient ring. 
We define an ideal $I_R$ of $R \otimes _ k R$,
\[
 I_R := ( x_i \otimes 1 -1 \otimes x_i ; i =1,\dots,r ) (R \otimes _k R).
\] 
Then we put 
\[
 P_R^n := R \otimes _k R / I_R^{n+1}
\]
and regard it as an $R$-module via the map 
\[
R \to  R \otimes _k R,\ a \mapsto a \otimes 1.  
\]
Then the module is finitely generated over $R$.
If $X := \Spec R$, then the $\cO_X$-module $\cP_X^n$ defined above is identified with 
the sheaf $\widetilde{P_R^n}$ associated to the $R$-module $P_R^n$.

Let $k[[\xx]] := k[[x_1,\dots ,x_r]]$ be a formal power series ring with $r$ variables
and $S := k[[\xx]]/\fb$ its quotient ring. Similarly we define an ideal $\hat I_S$ of $S \hat 
\otimes _ k S$,
\[
 \hat I_S := ( x_i \otimes 1 -1 \otimes x_i ; i =1,\dots,r ) (S  \hat \otimes _k S).
\] 
Then we put 
\[
 \hat P_S^n := S \hat \otimes _k S / \hat I_S^{n+1}
\]
and regard it as an $S$-module via the map 
\[
S \to  S \hat \otimes _k S,\ a \mapsto a \otimes 1.  
\]
Then the module is finitely generated over $S$.
For a scheme $Y = \Spec S$, we define a coherent $\cO _Y $-module $\hat \cP_Y ^n$ to be
the sheaf $ \widetilde{\hat P_S^n} $ associated to $\hat P_S^n$. 

\begin{defn}
Suppose that $Y$ is reduced and of pure dimension $d$, and that 
$\hat \cP_Y^n$ is locally free of constant rank $\binom{n+d}{d}$
on an open dense subset of $Y$.
Then we define the {\em $n$-th Nash blowup} of $Y$, denoted $\Nash _n (Y)$, to be 
$\Nash (Y , \hat \cP_Y ^n)$.
\end{defn}

Let $\hat \Delta^{(n)}_{Y} := \Spec \hat P^{n}_{S}$.
Then $\Nash _{n} (X)$ is identified with the union of 
the irreducible components of $\Hilb_{\binom{d+n}{d}}(\pr_{1}: \hat \Delta^{(n)}_{Y} \to Y)$
that dominate irreducible components of $Y$.

The condition that $\hat \cP_Y^n$ is locally free of constant rank $\binom{n+d}{d}$
on an open dense subset is probably superfluous. 
From the following lemma, when $Y$ is the completion of a variety at a point
or its irreducible component, the condition is, in fact, satisfied.

\begin{lem}\label{lem-compatible-Pn}
Let $R = k[\xx] /\fa$ and $\hat R := k[[\xx]] / \fa k[[\xx]]$.
Then there exists a natural isomorphism 
\[
    \hat P_ {\hat R}^n \cong P_R^n \otimes _R \hat R .
\]
\end{lem}

\begin{proof}
Let us view $k[\xx] \otimes _k k[\xx]$ (resp.\ $k[[\xx]] \hat \otimes _k k[[\xx]]$) 
as a $k[\xx]$-algebra (resp.\ a $k[[\xx]]$-algebra) by the map
$ x \mapsto x \otimes 1 $.
We have an isomorphism of $k[\xx]$-algebras,
\begin{align*}
\phi: k[\xx ] \otimes_k k[ \xx]  & \to k[\xx, \yy] :=k[\xx, y_1 , \dots, y_r] \\
          1 \otimes x_i & \mapsto x_i - y_i
\end{align*}
and an isomorphism of $k[[\xx]]$-algebras,
\begin{align*}
\hat \phi: k[[\xx ]] \hat \otimes _k k[[\xx]]  & \to k[[\xx, \yy]] :=k[[\xx, y_1 , \dots, y_r]] \\
         1 \otimes  x_i & \mapsto x_i - y_i.
\end{align*}
Then $ R \otimes _k R \cong  k[\xx, \yy] / \phi (\fa \otimes_k k[\xx] + k[\xx] \otimes_k \fa ) $ and 
\[
 P_{ R}^n \cong k[\xx,\yy]/ (\phi (\fa \otimes_k k[\xx] + k[\xx] \otimes_k \fa ) + (y_1,\dots,y_r)^{n+1}).
\]
Similarly, if $\hat \fa $ denotes $\fa k[[\xx]]$, then
\[
 \hat P_{ \hat R}^n \cong k[[\xx,\yy]]/ (\hat \phi (\hat \fa \hat \otimes_k k[[\xx]] + k[[\xx]] \hat 
 \otimes_k \hat \fa ) + (y_1,\dots,y_r)^{n+1}).
\]
We have  
\begin{align*}
&\hat R \otimes _{R} P_{R}^n \\
&\cong k[[\xx]][\yy] / ( \phi (\fa \otimes_k k[\xx] + k[\xx] \otimes_k \fa ) + (y_1,\dots,y_r)^{n+1}) \\
& \cong k[[\xx, \yy]] / (\hat \phi (\hat \fa \hat \otimes_k k[[\xx]] + k[[\xx]] \hat 
 \otimes_k \hat \fa )  + (y_1,\dots,y_r)^{n+1}) \\
&\cong \hat P_{\hat R}^n .
\end{align*}
\end{proof}

\begin{cor}\label{cor-formal-localization}
Let $X$ be a variety, $x \in X$ and $\hat X := \Spec \hat \cO_{X,x}$.
Then there exists a natural isomorphism
\[
 \Nash _n (\hat X) \cong \Nash _n (X) \times _X \hat X.
\]
\end{cor}

\begin{proof}
Let $f:\hat X \to X $ be a natural morphism.
From Lemma \ref{lem-compatible-Pn}, $\hat \cP_{\hat X}^n \cong f^* \cP_{ X}^n $, which
 implies the corollary.
\end{proof}

\section{General properties}\label{section-general-properties}

\subsection{Compatibility with etale morphisms}

\begin{thm}\label{thm-compatibility}
 Let $Y \to X$ be an etale morphism of varieties.
Then for every $n$, there exists a canonical isomorphism
\[
  \Nash_n (Y) \cong \Nash _n(X) \times _X Y . 
\]
\end{thm}

\begin{proof}
 Let $\Delta_X$ and $\Delta _Y $ be the diagonals in $X \times _k X$ and $Y \times _k Y$
 respectively.
Then the natural morphism
\[
 \Delta^{(n)}_Y \to \Delta _X^{(n)} \times _X Y
\]
is an isomorphism. This induces an isomorphism 
\[
 \Hilb _{\binom{d+n}{n}}(\pr_{1}:\Delta^{(n)}_Y \to Y ) 
 \cong  \Hilb _{\binom{d+n}{n}}(\pr_{1}:\Delta^{(n)}_X \to X ) \times _{X} Y 
\]
and the isomorphism of the assertion. 
\end{proof}

\subsection{Group actions}

Let $X$ be a variety of dimension $d$, $G$  an algebraic group over $k$ acting on $X$.
For each $l \in \NN$, we have a natural action of $G$ on $X \times _{k} \Hilb _l (X)$,
\begin{align*}
G \times _k X \times _{k} \Hilb _l  (X)  & \to X\times _{k} \Hilb_l (X) \\
(g ,x, [Z]) & \mapsto (gx,[g Z]) . 
\end{align*} 
When $l = \binom{d+n}{n}$, the subscheme $\Nash _{n}(X) \subseteq X \times _{k} \Hilb_{\binom{d+n}{n}}(X)$
is stable under this action. 
Thus the $G$-action on $X$ naturally lifts to $\Nash _n (X)$
and the morphism $\pi _n :\Nash _n (X) \to X$ is $G$-equivariant.

\subsection{Conductor and Jacobian ideals}

We now recall the conductor and Jacobian ideals, and their relation.
The conductor ideal plays an important role in what follows,
while the Jacobian ideal appears in Conjecture \ref{Yasuda-conjecture}.

Let $R$ be either a finitely generated $k$-algebra
or a local complete Noetherian ring with coefficient field $k$.
Suppose that $R$ is reduced and of pure dimension $d$.
Let $\tilde R$ be the integral closure of $R$ in the total ring of
fractions. 

\begin{defn}
The {\em conductor ideal} of $R$, denoted $\fc_{R}$,
is the annihilator of an $R$-module $\tilde R /R$.
\end{defn}

The conductor ideal is characterized as the largest ideal
of $R$ that is also an ideal of $\tilde R$.

\begin{defn}
When $R$ is a finitely generated $k$-algebra 
(resp.\ a complete local Noetherian ring with coefficient field $k$),
then the {\em Jacobian ideal} of $R$, denoted $\fj_{R}$, is
the $d$-th Fitting ideal of the module of K\"{a}hler differentials
 $\Omega_{R/k}$ (resp.\  the complete 
 module of K\"{a}hler differentials $\hat \Omega _{R/k}$). 
\end{defn}

If $R$ is represented as 
\[
R =k[x_{1},\dots,x_{m}] /(f_{1},\dots,f_{r})
\text{ or } R =k[[x_{1},\dots,x_{m}]] /(f_{1},\dots,f_{r}),
\]
then  $\fj_{R}$ is generated by the $(m-d)\times (m-d)$-minors
of the Jacobian matrix $ (\partial f_{i} / \partial x_{j} )_{i,j} $.

The conductor and Jacobian ideals commute with localizations.
Therefore they defines ideal sheaves on varieties. More directly,
if $X$ is a variety of dimension $d$ and $\nu :\tilde X \to X$ is the normalization, 
then the {\em conductor ideal sheaf} $\fc_{X} \subseteq \cO_{X}$ 
is defined to be the annihilator ideal sheaf of a coherent $\cO_{X}$-module $\nu_{*}\cO_{\tilde X}/\cO_{X}$.
The {\em Jacobian ideal sheaf} $\fj_{X} \subseteq \cO_{X}$ is defined to be 
the $d$-th Fitting ideal sheaf of the sheaf of K\"{a}hler differentials $\Omega_{X/k}$.
We call the closed subscheme $C_{X} \subseteq X$
defined by $\fc_{X}$ the {\em conductor subscheme}
and the closed subscheme $J_{X} \subseteq X$ defined by $\fj_{X}$ the {\em Jacobian subscheme}.
Similarly, when $X=\Spec R$ with $R$ a complete local Noetherian ring with
coefficient field $k$, then the conductor subscheme $C_{X} \subseteq X$
and the Jacobian subscheme $J_{X} \subseteq X$ are defined to be the subschemes
defined by $\fc_{R}$ and $\fj_{R}$ respectively.

The conducor and Jacobian ideals commute also with completion: 
Let $R$ is a finitely generated
$k$-algebra, $\fm \subseteq R$ is a maximal ideal,
and $\hat R$ the $\fm$-adic completion of $R$.
Then $\fj_{\hat R} = \fj_{R} \hat R$ and $\fc_{\hat R} = \fc_{R}\hat R$.

The relation of the conductor and Jacobian ideals is as follows:

\begin{thm}\label{thm-conductor-jacobian}
Let $R$ be either a finitely generated $k$-algebra or a local complete Noetherian 
ring with coefficient field $k$. Suppose that $R$ is reduced and of pure dimension $d$.
Then $\fj_{R} \subseteq \fc_{R}$.
\end{thm}

\begin{proof}
We prove only the case where $R$ is a finitely generated $k$-algebra.
The proof of the other case is parallel.

From the Noether normalization theorem,
there exists a $k$-homomorphism $\phi:k[x_{1},\dots,x_{d}] \to R$ 
which makes $R$ generically \'{e}tale over $k[x_{1},\dots,x_{d}]$. We can represent $R$ as 
\[
 R= k[x_{1},\dots,x_{d}] [x_{d+1},\dots,x_{m}]/(f_{1},\dots,f_{r}) = k[x_{1},\dots,x_{m}]/(f_{1},\dots,f_{r}).
\]
Then Lipman-Sathaye theorem \cite[Th.\ 2]{Lipman-Sathaye}
implies that 
every $(m-d) \times (m-d) $-minor of the matrix $(\partial f_{i}/\partial x_{j})_{\substack{1 \le i \le r \\ d+1 \le j \le m}}$
is contained in $\fc_{R}$. (For the case where $R$ is not a domain, see \cite[Th.\ 3.1]{Hochster}.)

For a suitable choice of variables $x_{1}, \dots, x_{m}$
and for every subset $\{j_{1},\dots,j_{d} \} \subseteq \{1,\dots,m\}$
of $d$ elements, $R$ is generically \'{e}tale over $k[x_{j_{1}},\dots,x_{j_{d}}]$.
Then $\fc_{R}$ contains every $(m-d) \times (m-d) $-minor
of the matrix $(\partial f_{i}/\partial x_{j})_{\substack{1 \le i \le r \\ 1 \le j \le m}}$.
Namely $\fj_{R} \subseteq \fc _{R}$.
\end{proof}

The following proposition is required in the following subsection:

\begin{prop}\label{prop-conductor-intersection}
Let $R$ be as above, $X := \Spec R$, and $X_{1},\dots,X_{l}$ be the irreducible components of $X$.
\begin{enumerate}
\item For $1 \le l' \le l$ and for $n \in \No$, if we put $X':= X_{1} \cup \dots \cup X_{l'}$, we have
\begin{align*}
C_{X'} \subseteq C_{X} \cap  X'  \text{ and} \\
J^{(n)}_{X'} \subseteq J^{(n)}_{X} \cap X' .
\end{align*}
Here $\subseteq, \cap, \cup$ are all scheme-theoretic.
\item 
The following inclusions hold
\[
 J_{X}  \supseteq C_{X}  \supseteq X_{1} \cap X_{2} .
\]
\end{enumerate}
\end{prop}

\begin{proof}
1. The first inclusion follows from the inclusion
\[
\widetilde{ X'} \subseteq \tilde X \times _{X} X'.
\]
To show the second one, it suffices to show $ J_{X'} \subseteq J_{X} \cap X' $.
If $R$ is finitely generated over $k$ and represented as
\[
R = k[x_{1},\dots,x_{m}] /(f_{1},\dots,f_{r}),
\]
and if 
\[
 R' =  k[x_{1},\dots,x_{m}] /(f_{1},\dots,f_{r},f_{r+1},\dots,f_{r'})
\]
is the coefficient ring of $X'$,
then $\fj_{R} R'$ is generated by 
the $(m-d) \times (m-d)$-minors of $(\partial f_{i} /\partial x_{j})_{\substack{1 \le i \le r \\ 1 \le j \le m}}$,
while $\fj_{R'}$ generated by the
$(m-d) \times (m-d)$-minors of $(\partial f_{i} /\partial x_{j})_{\substack{1 \le i \le r' \\ 1 \le j \le m}}$.
This shows the second inclusion of the assertion in this case.
The formal complete case is parallel. 

2.  The inclusion $J_{X} \supseteq C_{X}$ is equivalent to
Theorem \ref{thm-conductor-jacobian}. Concerning the other inclusion, 
from 1, we may suppose that $X_{1}$ and $X_{2}$ are the only 
irreducible components of $X$.
Let $R_{i}= R/I_{i}$, $i=1,2$, be the coefficient rings of
$X_{1}$ and $X_{2}$ respectively.
Since $R \subseteq R_{1} \times R_{2} \subseteq \tilde R$,
we have 
\[
\fc _{R} \subseteq \ann ( R_{1} \times R_{2} /R) = I_{1} + I_{2}.
\]
This prove the assertion.
\end{proof}

\subsection{Separation of analytic branches}

Let $X$ be a variety of dimension $d>0$, 
  $\hat X := \Spec \hat \cO_{X,x}$ the completion of $X$
at a point $x \in X$, and $\hat X_i$, $i=1,\dots,l$ its irreducible components.
Then we have 
\[
 \Nash _n (X) \times _X \hat X \cong \Nash _n (\hat X) \cong \bigcup _{i=1}^l \Nash _n (\hat X_i).
\] 
Let $[Z] \in \Nash _n (X)$ with $\pi_{n}([Z]) =x$. 
Then we can regard $[Z]$ also as a ($k$-)point of $\Nash _n (\hat X)$ and of $\Nash _n (\hat X_{i_0})$
for some $ 0 \le i_0 \le l$.  
Then  $Z$ is a closed subscheme of $\hat X_{i_0}$.
Moreover, from Corollary \ref{cor-contained-in-nth-nbhd}, $Z \subseteq \hat X_{i_0} \cap x^{(n)}$.

\begin{prop}\label{prop-general-separation}
Let $[Z] \in \Nash_{n} (X)$ with support $x$ and $Z \nsubseteq C_{X}$.
Then $Z$ is contained in a unique analytic branch  $\hat X_{i}$. 
Equivalently $ [Z] $ is contained in $\Nash _{n} (\hat X_{i})$ for a unique $i$. 
\end{prop}

\begin{proof}
From Proposition \ref{prop-conductor-intersection}, $Z$ can not be contained simultaneously 
in two irreducible components. This proves the proposition.
\end{proof}

Let $X,x,Z$ be as above. If $Z \nsubseteq J_{X}^{(d-1)}$, then from Theorem
\ref{thm-conductor-jacobian}, $Z \nsubseteq C_{X}$ and so $[Z] \in \Nash _{n} (\hat X_{i})$
for a unique $i$, say $i_{0}$.
Moreover from Proposition \ref{prop-conductor-intersection}, $Z \nsubseteq J^{(d-1)}_{\hat X_{i_{0}}}$.
As a consequence, Conjecture \ref{Yasuda-conjecture} is reduced to the following conjecture:

\begin{conj}
Let $R$ be a local complete Noetherian domain with coefficient field $k$
and $X := \Spec R$. Then for every $n \in \No$, 
the $n$-th Nash blowup $\Nash _{n} (X)$ is well-defined even if 
$X$ is not algebraizable. Moreover if 
$[Z] \in \Nash _{n} (X)$ with $Z \nsubseteq J_{X}^{(d-1)}$,  
then $\Nash _{n} (X)$ is regular at $[Z]$.
\end{conj}

The common idea in Proposition \ref{prop-general-separation}
and Conjecture \ref{Yasuda-conjecture} is that 
if $Z$ is too fat to be contained in some subscheme of $X$ 
like $C_{X}$ or $J_{X}^{(d-1)}$, then $X$ have mild singularities. 
The following lemma assures that a condition like 
$Z \nsubseteq C_{X}$ or $Z \nsubseteq  J_{X}^{(d-1)}$ 
holds for all $[Z] \in \Nash _{n} (X)$ if $n$ is sufficiently large,
and that $\Nash _{n} (X)$ has mild singularities everywhere.

\begin{prop}\label{prop-Hilbert-func-inclusion}
Let $X$ be a variety of dimension $d$ and
 $A \subseteq X$ a closed subscheme of dimension $<d$. 
Then there exists $n_0 \in \No$ such that for every $n \ge n_0$
and for every $[Z] \in \Nash_n (X) $, $Z \nsubseteq A$.
\end{prop}

\begin{proof}
Since $A$ is of dimension $ < d$, 
for every $a \in A $, the Hilbert function of $\cO_{A,a}$ is a polynomial of degree $<d$ for $n \gg 0$.
It follows that for $n \gg 0$, 
 \[
 \length \cO _{A,a}/ \fm _{A,a}^{n+1} < \binom{n+d}{d}.
 \] 
Because of the semi-continuity of Hilbert functions proved by Bennett \cite{Bennett},
for $n \gg 0$, the inequality holds simultaneously for all $a \in A$.

Let $[Z] \in \Nash _ n(X)$ and  $a$ its support. 
From Corollary \ref{cor-contained-in-nth-nbhd}, $Z \subseteq a^{(n)}$.
Since $\length \cO_Z = \binom{n+d}{d}$, if the inequality  holds,
then $Z \nsubseteq A \cap a^{(n)}$ and hence $Z \nsubseteq A$.
\end{proof}

\section{Higher Nash blowups of curves}\label{section-curves}

\subsection{A deformation-theoretic criterion for the regularity}

Let $R \subseteq k[[x]]$ be a complete $k$-subalgebra such that 
$k[[x]]$ is a finite $R$-module, $X := \Spec R$
and $\nu:\tilde X \to X$ its normalization. 
Since $X$ is algebraizable,  we can define higher Nash blowups of $X$.
To make computations below simpler, we fix the identification 
\[
 \tilde X = \Spec k[[y]]
\]
such that the ring homomorphism $ \nu^{*}: R \to k[[y]]$ is 
the composite of the inclusion $R \hookrightarrow k[[x]]$ and
the map $k[[x]] \to k[[[y]]] $, $x \mapsto -y$.
Then the complete fiber product of $X$ and $\tilde X$ is 
represented as
\[
 X \hat \times _{k} \tilde X: = \Spec R \hat \otimes _{k} k[[y]] = \Spec R[[y]].
\]
The graph $\Gamma_{\nu} \subseteq  X \hat \times _{k} \tilde X $ of $\nu$ is 
generically defined by $(x+y)$.
In precise, if $I \subseteq R[[y]]$ is the defining ideal of $\Gamma_{\nu}$,
then 
\[
 IR[[y]]_{I} = (x+y).
\]
Here $R[[y]]_{I}$ is the localization of $R[[y]]$ with respect to the prime ideal $I$.
Let $\cZ_{n} \subseteq X \hat \times _{k} \tilde X$ be the closed subscheme defined by
the $(n+1)$-th symbolic power of $I$,
\[
I^{(n+1)}:= R[[y]] \cap I ^{n+1} R[[y]]_{I}.
\]
Since the projection
\[
q_{n}:\cZ_{n} \to X
\]
is flat, 
 we obtain a corresponding birational morphism
\[
\phi_{n }:\tilde X \to \Nash_{n}(X)
\]
such that  $ \pi_{n} \circ \phi_{n} = \nu $.

Let $ o \in \tilde X $ be the closed point and $Z_n := q _n^{-1}(o) \subseteq X$,
the subscheme corresponding to $\phi _n (o) \in  \Nash _n (X)$.
Consider a natural morphism
\[
\epsilon : \Spec k[y] / (y^2) \to \Spec k[[y]] = \tilde X,
\]
which is a nonzero tangent vector of $\tilde X$ at $o$.
The fiber product 
\[
 \cZ_{n,\epsilon} :=  \cZ_n \times_ {q _n, \tilde X,\epsilon} \Spec k[y]/(y ^2)
 \subseteq X \times _k \Spec k[y]/(y ^2)
\]   is 
the first order embedded deformation of $Z_n \subseteq X$ corresponding to 
\[
 \phi_n \circ \epsilon :  \Spec k[y] / (y^2) \to \Nash _n (X).
 \]
Let $\fa_{n} \subseteq R$ be the defining ideal of $Z_{n}$, which is
identical to $ I^{(n+1)} $ modulo $(y)$.

\begin{thm}\label{thm-criterion-for-smoothness}
Suppose that $X$ is not regular.
Then the following are equivalent:
\begin{enumerate}
\item $\Nash _{n} (X)$ is regular.
\item $\phi_{n}$ is an isomorphism.
\item $\cZ_{n,\epsilon}$ is not the trivial embedded deformation of $Z_{n}$.
\item There exists an element $g \in I^{(n+1)} \subseteq R[[y]]$ such that
if we write 
\[
 g = g_{0} + g_{1} y + g_{2} y^{2} + \cdots , \ g_{i} \in R,
\]
then $g_{1} \notin \fa_{n}$.
\end{enumerate}
\end{thm}

\begin{proof}
$1 \Leftrightarrow 2$. Obvious.

$2 \Leftrightarrow 3$. 
The morphism $\phi_{n} \circ \epsilon $ corresponding to the 
pair $ ( \nu \circ \epsilon , \cZ_{n,\epsilon} ) $.
From the assumption, $\nu \circ \epsilon$ is the zero tangent vector, that is, factors as
$\Spec k[y]/(y^{2}) \to \Spec k \to X$.
Hence, $\phi_{n} \circ \epsilon$ is the zero tangent vector if and only if 
$\cZ_{n,\epsilon}$ is trivial.
This shows the equivalence $2 \Leftrightarrow 3$.

$3 \Leftrightarrow 4$.
If the defininig ideal of $\cZ_{n,\epsilon}$ in $R[y]/(y^{2})$ is generated by 
\[
 g_{j0} + g_{j1} y, \ g_{j0}, g_{j1} \in R , \ j= 1, \dots, m,
\]
then $\cZ _{n,\epsilon}$ corresponds to the homomorphism 
\begin{align*}
 \fa_n  \to R /\fa _n, \  g_{j0}  \mapsto g_{j1} .
\end{align*}
In particular, $\cZ_{n,\epsilon}$ being trivial is equivalent to
that the homomorphism is the zero map. Hence $3 \Leftrightarrow 4$.
\end{proof}

\begin{rem}
In the theorem above, the assumption that $X$ is not regular is unevitable.
For example, 
in characteristic $p >0$, if $X$ is regular, then 
$\cZ_{pm-1,\epsilon}$, $m \in \No$, are trivial. 
\end{rem}

\subsection{Associated numerical monoids}

A {\em numerical monoid} is by definition a submonoid $S$ of 
the (additive) monoid $\No$ with $\sharp (\No \setminus S) < \infty$.

To $ R \subseteq k[[x]]$ as above, we associate a numerical monoid 
\[
 S := \{ i \in \No | \exists f \in R, \ \ord f = i\}=\{0=s_{-1}<s_{0}<s_{1}< \cdots \}.
\]

\begin{thm}\label{thm-curve1}
Let $X : = \Spec R$. Suppose that $k$ has characteristic 0.
Then $\Nash _n (X)$ is regular if and only if $ s_n -1 \in S $.
\end{thm}

\begin{lem}\label{lem-matrices-regular}
Let  $\ba := \{ a_1  < a_2  < \cdots < a_{e}\} \subseteq \NN$  and define a $(e\times e)$-matrix
\[
M (n; \ba):=
\begin{pmatrix}
\binom{n}{a_1} & \binom{n}{a_1-1}  & \cdots &\binom{n}{a_1 - e+1} \\
\binom{n}{a_2} & \binom{n}{a_2-1}  & \cdots &\binom{n}{a_2 - e+1}  \\
  \vdots & \vdots &  \ddots & \vdots   \\ 
\binom{n}{a_e} & \binom{n}{a_e-1}  & \cdots & \binom{n}{a_e - e+1} 
\end{pmatrix}
\]
with entries in an algebraically closed field $k$ of characteristic 0.
 Here $ \binom{a}{b}:= 0$ if either $b > a$ or $b < 0$.
Then
\[
 \det M (n; \ba)=  \frac{ \prod _{ i < j} (a_j - a_i) \prod_{i=1}^e  
 \bigl((n + e -i) (n+e -i -1)\cdots (n + e - a_i ) \bigr) }
 { \prod _{i=1}^e a_i !} .
\] 
In particular, if $n+e-a_e >0$, then $\det M (n; \ba) \neq 0$, and the matrix $M(n;\ba)$
is regular.
\end{lem}

\begin{proof}
This matrix appears also in  \cite[page 353]{Arbarello-Cornalba-}.

Without changing the determinant, we can replace the first column with the sum of
the first and the second, and the second with the sum of the second and the third, and so on.
The resulting matrix is 
\[
\begin{pmatrix}
\binom{n+1}{a_1} & \binom{n+1}{a_1-1}  & \cdots&\binom{n+1}{a_1 - e+2} &\binom{n}{a_1 - e+1} \\
\binom{n+1}{a_2} & \binom{n+1}{a_2-1}  & \cdots&\binom{n+1}{a_2 - e+2} &\binom{n}{a_2 - e+1}  \\
  \vdots & \vdots &  \ddots &\vdots & \vdots   \\ 
\binom{n+1}{a_e} & \binom{n+1}{a_e-1}  & \cdots&\binom{n+1}{a_e - e+2} & \binom{n}{a_e - e+1} 
\end{pmatrix}.
\]
Again  we replace the first column with the second and the third, and so on.
We obtain 
\[
\begin{pmatrix}
\binom{n+2}{a_1} & \binom{n+2}{a_1-1}  & \cdots &\binom{n+2}{a_1 - e+3} &\binom{n+1}{a_1 - e+2} &\binom{n}{a_1 - e+1} \\
\binom{n+2}{a_2} & \binom{n+2}{a_2-1}  & \cdots&\binom{n+2}{a_2 - e+3} &\binom{n+1}{a_2 - e+2} &\binom{n}{a_2 - e+1}  \\
  \vdots        & \vdots             &  \ddots  &\vdots    &\vdots             & \vdots   \\ 
\binom{n+2}{a_e} & \binom{n+2}{a_e-1}  & \cdots&\binom{n+2}{a_e - e+3}&\binom{n+1}{a_e - e+2} & \binom{n}{a_e - e+1} 
\end{pmatrix}.
\]
Repeating this, we finally arrive at
\[
\begin{pmatrix}
 \binom{n+e-1}{a_1} &  \binom{n+e-2}{a_1 -1} & \cdots & \binom{n+1}{a_1 -e+2} &  \binom{n}{a_1-e+1}  \\
 \binom{n+e-1}{a_2} &  \binom{n+e-2}{a_2 -1} & \cdots & \binom{n+1}{a_2 -e+2} &  \binom{n}{a_2-e+1} \\
\vdots & \vdots & \ddots & \vdots &\vdots \\
  \binom{n+e-1}{a_e} &  \binom{n+e-2}{a_e -1} & \cdots & \binom{n+1}{a_e -e+2} &  \binom{n}{a_e-e+1} \\
\end{pmatrix} .
\]
(Check that this transformation makes sense even if the matrix $M(n;\ba)$ contains zero entries.)
Then we have
\begin{align*}
&\det M (n; \ba) \\ 
&= \det  
\begin{pmatrix}
 \binom{n+e-1}{a_1} &  \binom{n+e-2}{a_1 -1} & \cdots & \binom{n+1}{a_1 -e+2} &  \binom{n}{a_1-e+1}  \\
 \binom{n+e-1}{a_2} &  \binom{n+e-2}{a_2 -1} & \cdots & \binom{n+1}{a_2 -e+2} &  \binom{n}{a_2-e+1} \\
\vdots & \vdots & \ddots & \vdots &\vdots \\
  \binom{n+e-1}{a_e} &  \binom{n+e-2}{a_e -1} & \cdots & \binom{n+1}{a_e -e+2} &  \binom{n}{a_e-e+1} \\
\end{pmatrix} \\
& =   \prod_{i=1}^e  
 \bigl((n + e -i) (n+e -i -1)\cdots (n + e - a_i ) \bigr) \times \\*
 & \qquad \det \begin{pmatrix}
 ( a_1 !)^{-1} &  ( (a_1-1) !)^{-1}  & \cdots  &  ( (a_1-e+1) !)^{-1}  \\
 ( a_2 !)^{-1} &  ( (a_2-1) !)^{-1}  & \cdots  &  ( (a_2-e+1) !)^{-1}\\
\vdots & \vdots & \ddots &  \vdots \\
 ( a_e !)^{-1} &  ( (a_e-1) !)^{-1}  & \cdots  &  ( (a_e-e+1) !)^{-1}\\
\end{pmatrix} \\
&= \frac{  \prod_{i=1}^e  
 \bigl((n + e -i) (n+e -i -1)\cdots (n + e - a_i ) \bigr) }
 { \prod _{i=1}^e a_i !} \times \\*
 & \qquad \det \begin{pmatrix}
 1 &  a_1  & a_1 (a_1 -1) &  \cdots & a_1 (a_1 -1) \cdots (a_1 -e + 2)  \\
 1 &  a_2  & a_2(a_2 -1) &  \cdots & a_2 (a_2 -1) \cdots (a_2 -e + 2)  \\
\vdots & \vdots & \vdots & \ddots &  \vdots \\
 1 &  a_e & a_e (a_e -1) &  \cdots & a_e (a_e -1) \cdots (a_e -e + 2)  \\
\end{pmatrix} \\
& =\frac{  \prod_{i=1}^e  
 \bigl((n + e -i) (n+e -i -1)\cdots (n + e - a_i ) \bigr) }
 { \prod _{i=1}^e a_i !} \times \\*
 & \qquad \det \begin{pmatrix}
 1 &  a_1  & a_1 ^2 &  \cdots & a_1^{e-1} \\
 1 &  a_2  & a_2 ^2 &  \cdots & a_2^{e-1}  \\
\vdots & \vdots & \vdots & \ddots &  \vdots \\
  1 &  a_e  & a_e ^2 &  \cdots & a_e^{e-1}   \\
\end{pmatrix} \\
&=\frac{  \prod_{i=1}^e  
 \bigl((n + e -i) (n+e -i -1)\cdots (n + e - a_i ) \bigr) }
 { \prod _{i=1}^e a_i !} \prod _{ i < j} (a_j - a_i).  \\
 &  \text{(Vandermonde's determinant)} 
\end{align*}
\end{proof}

\begin{proof}[Proof of Theorem \ref{thm-curve1}]
Put $ T := \No \setminus S = \{t_{1} < t_{2} < \cdots < t_{l} \} $.
Let $\bt _{n,0} := \{ t \in T | t  < s_{n}\}=\{t_{1}<t_{2}<\cdots < t_{l_{n}}\}$, where $l_{n} := \sharp \bt_{n,0}$,
 and $\bu_{n,0}:= \bt_{n,0} \cup \{s_{n}\}$.
Then $s_{n}=l_{n} +n +1$.
From Lemma \ref{lem-matrices-regular}, the matrix $M (n+1;\bu_{n,0})$ is regular.
We define $r_{n,i} \in k$, $i=1,\dots,l_n$, by the equation
\[ 
M (n+1; \bu_{n,0} )
\begin{pmatrix}
r_{n,0}\\
\vdots \\ 
r _{n,l_{n}-1} \\
r_ {n,l_n} \\
\end{pmatrix}
=  
\begin{pmatrix}
0 \\
\vdots \\
0 \\
1 \\
\end{pmatrix}.
\]
Then we define a homogeneous polynomial of degree $s_n$,
\[
f_{n,0}:= (r_{n,0} y ^ {l_n} + r_{n,1} x y^{l_n-1} + \cdots + r_{n,l_n-1}  x^{l_n-1} y + r_{n,l_n} x^{l_n}) 
( x + y)^{n+1} \in k[x,y].
\]
For $1 \le i \le l_n$,
the coefficient of $x^{t_i}y^{s_n-t_i}$ in $f_{n}$ is
\[
 r_{n,0} \binom{n+1}{t_i} + r_{n,1} \binom{n+1}{t_i -1} + \cdots + r_{t,l_n} \binom{n+1}{t_i -l_n} = 0 
\]
and the coefficient of $x^{s_n} = x^{n+l_n+1} $ is $1$.

For $ j \in \NN $, 
we put $\bt_{n,j}:=\{ t \in T | t \le s_{n} +j \}$. 
If $m_{n,j} := l_{n} +j -\sharp \bt_{n,j} \ge 0 $, 
then we put $\bs_{n,j}:= \{s_{0} < s_{1} < \cdots < s_{m_{n,j}}\}$
and $\bu _{n,j}:= \bt_{n,j} \cup \bs_{n,j}$. 
Then $\sharp \bu_{n,j} = l_{n} + j +1$. Since
\[
 (n+1) + \sharp \bu_{n,j} -\max \bu_{n,j} \ge (n+1) + l_{n}+j+1 -(s_{n}+j) =1,
\]
from Lemma \ref{lem-matrices-regular}, $M(n+1 , \bu_{n,j})$ is regular. 
Therefore from the same argument as above,
for every $(d_{i};i \in \bu _{n,j}) \in k^{\bu_{n,j}}$,
 there exists a unique homogeneous polynomial $g \in k[[x,y]]$ 
 such that 
\begin{enumerate}
\item $g$ has degree $(l +j )+(n+1) = s_{n}+j$,
\item $g$ is divided by $(x+y)^{n+1}$, and
\item for each $i \in \bu_{n,j}$, the coefficient of the term $x^{i}y^{s_{n}+j-i}$ is $d_{i}$.
\end{enumerate}

Now we inductively choose homogeneous polynomials 
$f_{n,j}$, $j \in \NN$, of degree $s_{n}+j$ divisible by $(x+y)^{n+1}$ as follows:
For each $i \in \No$, we can take an element
\[
h_{i} = \sum_{j \ge 0} h_{i,j} x^{i + j}  \in R ,\ h_{i,j} \in k
\]
such that 
\begin{enumerate}
\item $h_{0}=1$,
\item for $i \in S$, $h_{i,0}=1$, 
\item for $i \in T$, $h_{i}=0$, and 
\item if $j>0$ and $i+j \in S$, then $h_{i,j}=0$.
(In particular, if $i > t_{l}$ and $j >0$, then $h_{i,j}=0$.)
\end{enumerate}
Suppose that we have chosen $f_{n,0},f_{n,1},\dots,f_{n,j-1}$.
 Let
$c_{i,j'}$, $i \le s_{n} + j'$, $0 \le j' < j$, be the coefficient of $x^{i}y^{s_{n}+j'-i}$ in $f_{n,j'}$.
By convention, we put $c_{i,j'}:=0$ for $i<0$ or for $j' <0$.
For $i \in \bs_{n,j}$, put $c_{i,j} := 0$.
For $i \in \bt_{n,j} $, put 
\[
 c_{i,j} := \sum _{ a = 1}^{j} c_{i-a, j-a} h _{i-a,a}.
\]
Then we choose $f_{n,j}$ such that for every $i \in \bu_{n,j}$,
the coefficient of $x^{i}y^{s_{n}+j-i}$ is $c_{i,j}$.

We claim that for $j \gg 0$, $f_{n,j}=0$.
To see this, we first observe that for $j \gg 0$,
the coefficients of $x^{i}y^{s_{n}+j-i}$, $i \in \{s \in S | s < t_{l}\} \subseteq \bs_{n,j}$,
are all $0$.  Then if necessary, replacing $j$ with a still larger integer, 
we obtain that $f_{n,j-1},f_{n,j-2},\dots,f_{n,j-t_{l_{n}}}$ all
have this property. Then for every $i \in \bt_{n,j}$, $c_{i,j}=0$.
From the uniqueness, $f_{n,j}=0$.

Define  $f_{n} := \sum _{j=0}^{\infty} f_{n,j}$.  
Then $f_{n}$ is divided by $(x+y)^{n+1}$. 
Moreover,
\begin{align*}
f_{n} & =  \sum _{i,j} c_{i,j}x^{i} y ^{s_{n}+j-i} \\
 & =\sum_{\substack{i,j \\ i \in S}} c_{i,j} x^{i}y^{s_{n}+j-i}
 + \sum _{\substack{i,j \\ i \in T}}  (\sum _{ a = 1}^{j} c_{i-a, j-a} h _{i-a,a}) x^{i}y^{s_{n}+j-i} \\
 &= \sum_{\substack{i,j \\ i \in S}} c_{i,j} x^{i}y^{s_{n}+j-i}
 + \sum _{\substack{i,j,a \\ a>0, i+a \in T}}  c_{i, j} h _{i,a} x^{i+a}y^{s_{n}+j-i} \\
 &= \sum_{i,j} c_{i,j} h_{i} y^{s_{n}+j-i} .
\end{align*}
Thus $f_{n} \in R[[y]]$ and so $f_{n} \in I^{(n+1)}$.
By construction, 
\[
f_{n}(x,0) = x^{s_{n}} + (\text{higher terms}) \in \fa_{n}.
\]
Similarly for every $n' \ge n$, $f_{n'} \in I^{(n+1)}$,
and 
\[
f_{n'}(x,0) = x^{s_{n'}} + (\text{higher terms}) \in \fa_{n} .
\]
Since 
\[
\length R/ (f_{n'}(x,0);n' \ge n) = n+1,
\]
$\fa_{n}$ is in fact generated by $f_{n'}(x,0)$, $n' \ge n$, 
and identical to $\{f \in R | \ord f \ge s_{ n}\}$.
It follows that $I^{(n+1)}$ is generated by
$f_{n'}$, $n' \ge n$. 

Write 
\[
f_{n} \equiv f_{n}(x,0) + g_{n} y  \mod (y^{2}), \ g_{n} \in R.
\]
From Theorem \ref{thm-criterion-for-smoothness}, $\Nash_{n}(X)$ is regular
if and only if for some $n' \ge n$, $g_{n'}\notin \fa_{n}$.
For every $n' >n$, $g_{n'}$ has order $\ge s_{n}$,
and so $g_{n'} \in \fa_{n}$.

Now $\Nash_{n}(X)$ being regular
is equivalent to that $g_{n}$ has order $s_{n}-1$,
or equivalently $c_{s_{n}-1 ,0} \ne 0$.
If $s_{n}-1 \in T$, 
then $s_{n}-1 = t_{l_{n}}$ and by the construction,
$c_{s_{n}-1 ,0}=0$.
If  $s_n -1 \in S$, then put $ \bu' _{n,0} := \{ t_1, \dots , t_{l_n} , s_n -1 \} $.
From Lemma \ref{lem-matrices-regular}, the matrix $M(n+1;\bu ' _{n,0})$ is regular.
We have
\[
M(n+1;\bu'_{n,0})
 \begin{pmatrix}
r_{n,0}\\
\vdots \\ 
r _{n,l_n-1} \\
r_ {n,l_n} \\
\end{pmatrix}
=
 \begin{pmatrix}
0\\
\vdots \\ 
0 \\
c_{s_{n}-1 ,0}\\
\end{pmatrix} \ne 0 .
\]
This completes the proof.
\end{proof}

\begin{expl}\label{expl-5,7}
We note that for every numerical monoid $S$,
there exists $R \subseteq k[[x]]$ whose associated monoid is $S$.
Suppose that $ S $ is the numerical monoid generated by $5$ and $7$.
Then
\[
 S = \{0,5,7,10,12,14,15,17,19,20,21,22, n; n \ge 24 \}.
\]
Theorem \ref{thm-curve1} now says that
 \[
 \Nash _n (X) \text{ is }
 \begin{cases}
\text{singular} & (n=0,1,2,3,4,6,7,11) \\
\text{regular} & \text{(otherwise)}.
\end{cases} 
 \]
\end{expl}

\begin{expl}
If for some $m$,  $S= \{0,  m,m+1,m+2 , \dots \}  $, 
then for every $n >0$, $\Nash _n (X)$ is regular. 
\end{expl}

\subsection{Conjecture \ref{Yasuda-conjecture} for curves}

\begin{cor}\label{cor-general-curve}
Suppose that $k$ has characteristic $0$.
Let $X$ be either a variety of dimension 1 or 
$\Spec R$ with $R$ a reduced local complete Noetherian ring with coefficient field $k$.
Let $C \subseteq X$ be the conductor subscheme
and $[Z] \in \Nash _{n} (X)$ with $Z \nsubseteq C$.
Then $\Nash_{n} (X)$ is regular at $[Z]$. 
In particular, Conjecture \ref{Yasuda-conjecture} is true
in dimension $1$.
\end{cor}

\begin{proof}
The second assertion is a consequence of the first and 
Theorem \ref{thm-conductor-jacobian}. We will now prove the
first assertion. 
From Corollary \ref{cor-formal-localization}, we may suppose that 
$X = \Spec R$ with $R$ a reduced local complete Noetherian ring with coefficient field $k$.
From Proposition \ref{prop-conductor-intersection}, 
$Z$ is contained in a unique irreducible component of $ X$,
say $X_{0}$. If $C_{0}$ is the conductor subscheme
of $ X_{0}$, then from Proposition \ref{prop-conductor-intersection},
we have $Z \nsubseteq C_{0}$.
Hence it suffices to prove only the case where $R$ is a domain,
the case as in Theorem \ref{thm-curve1}.
 With the notations as in Theorem \ref{thm-curve1},
 the conductor ideal $\fc$ of $R$ is $(x^{i};  i > t_{l})$. 
If $s_{n} > t_{l} +1$,
then from Theorem \ref{thm-curve1}, $\Nash_{n}(X)$ is smooth.
If $ s_{n} \le t_{l}+1 $, then as we saw in the proof of Theorem \ref{thm-curve1},
$\fa_{n} = \{ f \in R | \ord f \ge s_{n} \} \supseteq \fc $ and 
the condition, $Z \nsubseteq C$, is not satisfied. 
This completes the proof.
\end{proof}

If $Z = C$, then $\Nash _{n } (X)$ is not generally smooth at $[Z]$.
In fact, with $X=\Spec R$ as in Theorem \ref{thm-curve1}, 
if $n = i_{0} := \max \{i | s_{i}-1 \notin R \}$, then 
$Z_{n} =C$ and $\Nash _{n} (X)$ is not smooth at $[Z_{n}]$. 
Therefore we can not replace $C$ in the corollary with any smaller 
subscheme of $X$.

\subsection{Positive characteristic}

As the following propositions show, it is impossible to resolve curve
singularities in positive characteristic via higher Nash blowups.

\begin{prop}\label{prop-positive-char-curve}
Let $X = \Spec R$ be as in Theorem \ref{thm-curve1}.
Suppose that $k$ has  characteristic $p >0$.
Then for $e \gg 0$, 
\[
\Nash_{p^e-1} (X)\cong X.
\]
\end{prop}

\begin{proof}
For $e \gg 0$, 
\[
 (x + y)^{p^e} = x^{p^e}+y^{p^e} \in R \hat \otimes_k R \subseteq k[[x,y]].
\] 
Let $\cW \subseteq X \times _k X$ be the closed subscheme defined by the ideal 
$  ( x+y)^{p^e}$.
If $q \in X$ is the image of the origin $o \in \tilde X = \Spec k[x]$,
then the  fiber of $\pr _2 :\cW \to X$ over $q$ is $\Spec R/ x^{p^e}R$.
From \cite[Lem.\ 11.12]{Eisenbud}, 
\[
\length R/ x^{p^e}R = \length k[x]/x^{p^e}k[x] = p^e .
\]
From \cite[Ex.\ 20.13]{Eisenbud}, $\pr _1 :\cW \to X$ is flat.
There exists a corresponding morphism
\[
 X \to \Nash _{p^e-1} (X) ,
\] 
which is the inverse of $\pi_{p^e-1}: \Nash_{p^e-1}(X) \to X$.
We have proved the proposition.
\end{proof}

\begin{prop}\label{prop-cusp-char2}
Suppose that $k$ has characteristic either 2 or 3.
Let $X := \Spec k[[x^2,x^3]]$.
Then for every $n \in \No$, $\Nash _n (X) \cong X$.
\end{prop}

\begin{proof}
We first consider the case of characteristic 2.
For $ n \in \NN$,
\begin{align*}
 &(x+y)^{n} \\
 & = x ^{n} + n x^{n-1} y +  
( \sum_{i=2}^{n-2}\binom{n}{i}  x ^i y^{n-i} ) + nxy^{n-1}+y^{n} \\
 & = \begin{cases}
 x ^{n} +  
 ( \sum_{i=2}^{n-2}\binom{n}{i}  x ^i y^{n-i} )+y^{n} & (n \text{ even}) \\
 x ^{n} +  x^{n-1} y +  
 ( \sum_{i=2}^{n-2}\binom{n}{i}  x ^i y^{n-i} ) + xy^{n-1}+y^{n} & (n \text{ odd}).
\end{cases}
\end{align*}
Thus for odd $n$, $(x+y)^{n+1} \in R \hat \otimes _k R$. By the same argument with the proof of the last
proposition, we see that $\Nash _n (X) \cong X$.

For even $n$, 
the coefficients of $x ^{n+2}y$ and $xy^{n+2}$ in $(x^2 + xy +y^2) (x+y)^{n+1}$
are both zero. Therefore $ (x^2 + xy +y^2) (x+y)^{n+1} \in R \hat \otimes _k R$.
For an ideal
\[
 I := ( (x+y)^{n+2} , (x^2 + xy +y^2) (x+y)^{n+1}   ) \subseteq R \hat \otimes _k R ,
\]
we have 
\[
\length R / I R = \length R / (x^{n+2},x^{n+3}) R =  n+1.
\]
Again by the same argument,
 we can show the assertion in the case where $n$ is even.
 
We next consider the case of characteristic 3. 
Similarly we have
\begin{align*}
(x+y)^n &\in R \otimes _k R  \ (n \equiv 0  \mod 3), \\
(x-y)(x+y)^n = (x^2-y^2)(x+y)^{n-1} &\in R \otimes _k R  \ (n \equiv 1 \mod 3), \\
( x^2 +  xy +y^2 ) (x+y)^n & \in R \otimes _k R \  (n \equiv 2  \mod 3).
\end{align*}
For each $n \in \NN$, we define an ideal   $I \subseteq R \otimes _k R$ as follows:
\begin{align*}
I := \begin{cases}
((x-y)(x+y)^{n+1}, (x+y)^{n+3}  )&  (n \equiv 0  \mod 3) \\
((x+y)^{n+2}, ( x^2 +  xy +y^2 ) (x+y)^{n+1})& (n \equiv 1  \mod 3) \\
((x + y)^{n+1}) & (n \equiv 2  \mod 3).
\end{cases}
\end{align*}
Then 
\begin{align*}
\length R/IR  & = 
\begin{cases}
\length R /x^{n+1} R \text{ or}\\
\length R/(x^{n+2},x^{n+3}) R
\end{cases} \\
& = n+1.
\end{align*} 
We can similarly show the assertion.
\end{proof}

%%%%%%%%%%%%%%%%%%%%%%%%%%%%%%%%%%%%%%%%%%%%%%%%%%%%%%%%%%%%%%%%%%%%%%%%%%%%%%%

\end{document}